\documentclass[11pt]{article}
\usepackage{latexsym}
\usepackage{amssymb}
\usepackage{amsmath,amscd}
\usepackage{amsthm}
\usepackage{mathrsfs}
\allowdisplaybreaks

\newtheorem{theorem}{Theorem}[section]
\newtheorem{lemma}{Lemma}[section]
\newtheorem{proposition}{Proposition}[section]
\newtheorem{definition}{Definition}[section]
\newtheorem{remark}{Remark}[section]
\newtheorem{corollary}{Corollary}[section]
\newtheorem{example}{Example}[section]

\begin{document}

\title{\bf The random case of Conley's
theorem: III. Random semiflow case and\\ Morse
decomposition\thanks{Published in: {\em Nonlinearity} {\bf 20}
(2007), 2773--2791.}}

\author{Zhenxin Liu\\
{\small College of Mathematics, and Key Laboratory of Symbolic
Computation and Knowledge} \\ {\small  Engineering of Ministry of
Education, Jilin University,  Changchun 130012,
P.R. China}\\
{\small E-mail: zxliu@jlu.edu.cn}}

\date{}
\maketitle

\begin{abstract}
In the first part of this paper, we generalize the results of the
author \cite{Liu,Liu2} from the random flow case to the random
semiflow case, i.e. we obtain Conley decomposition theorem for
infinite dimensional random dynamical systems. In the second part,
by introducing the backward orbit for random semiflow, we are able
to decompose invariant random compact set (e.g. global random
attractor) into random Morse sets and connecting orbits between
them, which generalizes the Morse decomposition of invariant sets
originated from Conley \cite{Con} to the random semiflow setting and
gives the positive answer to an open
problem put forward by Caraballo and Langa \cite{CL}. 
\end{abstract}

\section{Introduction}

This paper is the third and also the final part of the series of
papers \cite{Liu,Liu2}, which aim at studying, in random setting,
Conley decomposition theorem and Morse decomposition theorem. Both
of these two theorems are originated from Conley \cite{Con}. My
other two related (joint) works are \cite{Liu-in,Liu3}.

Conley's fundamental theorem of dynamical systems \cite{Con} can be
stated as follows.

\begin{theorem}\label{Conley}{\bf (Conley's fundamental theorem of dynamical systems).}
Any flow on a compact metric space decomposes the space into a chain
recurrent part and a gradient-like part.
\end{theorem}

\noindent For the importance of the theorem, it was adapted for maps
on compact spaces by Franks \cite{Fra}, was later established for
maps on locally compact metric spaces by Hurley \cite{Hu0,Hu1}, and
was extended by Hurley \cite{Hu2} for semiflows and maps on
arbitrary metric spaces. Recently, the author \cite{Liu,Liu2}
extended Conley decomposition theorem to random dynamical systems
(RDS) on Polish spaces. The results of \cite{Liu,Liu2} was written
for finite dimensional (i.e. random flow) case, so a natural
question is whether the theorem also holds for infinite dimensional
RDS (i.e. random semiflow). This is just what we answer in the first
part of the present paper. Through simple observation, we can
characterize the random chain recurrent set in terms of random
attractors similar to \cite{Liu,Liu2}. But when we consider the
complete Lyapunov function for random semiflow, the construction in
\cite{Liu2} (see expression (4.1) in \cite{Liu2}) is not applicable.
The reason is as follows: 1) random semiflow is not defined for
$t<0$. 2) The random attractor is only forward invariant. Even if we
introduce backward orbits for random semiflow as we do in Section 4,
we can not still define $\tau(\omega,x)=-\infty$ when $x\in
A(\omega)$ (see (4.1) in \cite{Liu2}) because we can not conclude
that there must be a backward orbit through $x$ which lie on the
attractor. 3) Since backward orbits are not necessarily unique,
$\tau(\omega,x)$ is not well defined. To bypass the obstacles
mentioned above, we will construct a new complete Lyapunov function
following Conley \cite{Con} as well as Arnold and Schmalfuss
\cite{Ar2}, which has weaker properties than the complete Lyapunov
function in \cite{Liu2}.

Morse decomposition theorem, originated from Conley \cite{Con}, is
very useful in studying the inner structure of invariant sets, e.g.
global attractor (see {\cite{Hale,Tem} for comprehensive study of
it), which can be stated at the abstract level as follows.

\begin{theorem}{\bf (Morse decomposition theorem).}
Any flow restricted to an invariant compact set decomposes the
compact set into finite number of invariant compact subsets (i.e.
Morse sets) and connecting orbits between them.
\end{theorem}
\noindent To be more specific, let $\varphi$ be a flow and $S$ be an
invariant compact set of $\varphi$. Assume that
$(A_i,R_i),i=1,\cdots,n$ are attractor-repeller pairs of $\varphi$
with
\[
\emptyset=A_0\varsubsetneq A_1\varsubsetneq\cdots\varsubsetneq
A_n=S~{\rm and}~S=R_0\varsupsetneq
R_1\varsupsetneq\cdots\varsupsetneq R_n=\emptyset.
\]
Then the family $D=\{M_i\}_{i=1}^{n}$ of invariant compact subsets
of $S$, defined by
\[
M_i=A_i\cap R_{i-1},~1\le i\le n,
\]
is called a {\em Morse decomposition} of $S$, and each $M_i$ is
called {\em Morse set}. If $D$ is a Morse decomposition, $M_D$ is
defined to be $\bigcup_{i=1}^nM_i$, which completely describes the
asymptotic behaviors of $\varphi$. For detailed definitions and
further properties, see \cite{Con}.

For example, consider the differential equation $\dot{x}=(1-x^2)x$
with $x\in\mathbb R$ and assume that $\varphi$ is the flow generated
by it, see  Figure \ref{fig}. It is clear that $\mathcal A:=[-1,1]$
is the global attractor of $\varphi$. By restricting $\varphi$ to
$\mathcal A$, it is easy to see that all attractors in $\mathcal A$
are $\emptyset$, $\mathcal A$, $\{-1\}$, $\{1\}$, $\{-1,1\}$. If we
set $A_0=\emptyset$, $A_1=\{-1\}$, $A_2=\{-1,1\}$, $A_3=\mathcal A$,
then the corresponding repellers are $R_0=\mathcal A$, $R_1=[0,1]$,
$R_2=\{0\}$, $R_3=\emptyset$ and hence the corresponding Morse sets
are $M_1=\{-1\}$, $M_2=\{1\}$, $M_3=\{0\}$. Therefore,
$D=\{M_1,M_2,M_3\}$ is a Morse decomposition of $\mathcal A$.

\begin{figure}\label{fig}
\setlength{\unitlength}{0.1in} \centering
\begin{picture}(40,10)
\put(0,5){\vector(1,0){5}} \put(5,5){\line(1,0){10}}
\put(10,5){\circle*{0.5}} \put(20,5){\vector(-1,0){5}}
\put(20,5){\vector(1,0){5}} \put(25,5){\line(1,0){10}}
\put(40,5){\vector(-1,0){5}} \put(20,5){\circle*{0.5}}
\put(24,5){\line(1,0){6}} \put(30,5){\line(-1,0){6}}
\put(30,5){\circle*{0.5}} \put(10,4){\makebox(0,0){$-1$}}
\put(10,6){\makebox(0,0){$M_1$}} \put(20,4){\makebox(0,0){$0$}}
\put(20,6){\makebox(0,0){$M_3$}} \put(30,4){\makebox(0,0){$1$}}
\put(30,6){\makebox(0,0){$M_2$}}
\end{picture}\caption{A Morse decomposition for the global attractor of $\dot{x}=(1-x^2)x$}
\end{figure}
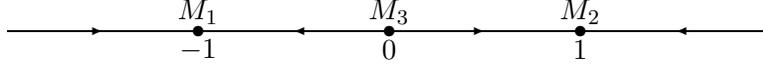

Similar to the deterministic case, the global random attractor is a
very crude object as it is the largest invariant random compact set,
which means that it includes all the smaller invariant random sets
of the systems and it also contains all the unstable manifolds
associated to all these invariant parts of the attractor. The
dynamics on the global random attractor is generally so complicated
that the structure of global random attractor is really important
for us to understand the asymptotic behavior of the system. In fact,
this theory, closely related to stochastic bifurcation theory, ``is
still in its infancy" (Arnold \cite{Ar1}). Some specific models have
been studied in this aspect, see, for example, Arnold and Boxler
\cite{AB}, Baxendale \cite{Ba}, Keller and Ochs \cite{KO},
Schenk-Hopp\'e \cite{Sche} etc for the finite dimensional case.
Caraballo et al \cite{CLR} studied the first stochastic bifurcation
problem in an infinite dimensional case, i.e. they showed that the
global random attractor of Chafee-Infante equation with a
multiplicative noise undergoes a stochastic pitchfork bifurcation
when the coefficient of linear term passes through the first
eigenvalue of the negative Laplacian from below. For general
theoretical results, Ochs \cite{Och} firstly considered the Morse
decomposition of weak random attractors; Crauel et al \cite{Cra} and
Liu et al \cite{Liu3} considered the Morse decomposition for random
flow on compact metric space (it is clear that their results also
hold when we consider Morse decomposition of global random
attractors instead of the entire state space). But these results are
written only in random flow case, i.e. finite dimensional case; ``it
is at all not clear if something similar (i.e. Morse decomposition
for global random attractors) could be true for SPDEs" (Caraballo
and Langa \cite{CL}). The obstacles to this open problem may contain
the following: 1) random backward orbits (see Section 4 for the
definition) are not introduced for random semiflows to our best
knowledge. 2) Generically, there is no uniqueness for random
backward orbits. 3) Even for a point (or random variable) in an
invariant random set, there may be backward orbit which does not lie
on the invariant set through the point (or random variable). In the
second part of the paper, by introducing random backward orbits and
careful treating the difference between random semiflow and random
flow, we give the positive answer to this open problem, i.e. we
obtain that Morse decomposition for global random attractors in
infinite dimensional case also holds. In fact, we prove that the
corresponding result holds for any invariant random compact set.

\section{Preliminaries}

Throughout the paper, we assume that $X$ is a Polish space, i.e. a
separable complete metric space. In this section, we will give some
preliminaries for later use.  Firstly we give the definition of
continuous random dynamical systems (cf Arnold \cite{Ar1}).

\begin{definition}
Let $X$ be a metric space with a metric $d_X$. A {\em random
dynamical system (RDS)}, shortly denoted by
$\varphi$, consists of two ingredients: \\
(i) A model of the noise, namely a metric dynamical system $(\Omega,
\mathcal F, \mathbb P, (\theta_t)_{t\in \mathbb R})$, where
$(\Omega, \mathcal F, \mathbb P)$ is a probability space and $(t,
\omega)\mapsto \theta_t\omega$ is a measurable flow which leaves
$\mathbb P$ invariant, i.e. $\theta_t\mathbb P=\mathbb P$ for all
$t\in \mathbb R$.\\
(ii) A model of the system perturbed by noise, namely a cocycle
$\varphi$ over $\theta$, i.e. a measurable mapping $\varphi: \mathbb
R^+\times \Omega\times X \rightarrow X,
(t,\omega,x)\mapsto\varphi(t,\omega,x)$, such that
$\varphi(t,\omega,\cdot)=\varphi(t,\omega):X\rightarrow X$ is
continuous for arbitrary $t\ge 0$ and $\omega\in\Omega$, moreover
\begin{equation}\label{phi}
\varphi(0,\omega)={\rm id}_X, \varphi(t+s,\omega)=\varphi(t,\theta_s
\omega)\circ\varphi(s,\omega)\quad {\rm for ~all}\quad t,s\in\mathbb
R^+,\omega\in\Omega.
\end{equation}
\end{definition}

In \cite{Liu,Liu2}, $\varphi$ is defined for every $t\in\mathbb R$
(i.e. random flow), which is usually true for the RDS generated by
finite dimensional random and stochastic differential equations,
i.e. random and stochastic ODEs. Here, $\varphi$ is only defined for
$t\ge0$ (i.e. random semiflow), which is usually generated by random
and stochastic PDEs.

\begin{definition}\label{invariant}
A random set $D$ is said to be {\em forward invariant} under the RDS
$\varphi$ if $\varphi(t,\omega)D(\omega)\subset D(\theta_t\omega)$
for all $t\ge 0$ almost surely; It is said to be {\em invariant} if
$\varphi(t,\omega)D(\omega)=D(\theta_t\omega)$ for all $t\ge 0$
almost surely.
\end{definition}

\begin{definition}
Assume that $D$ is a random set, then the {\em omega-limit set of
$D$}, $\Omega_D$, is defined to be
\[
\Omega_D(\omega):=\bigcap_{t\ge 0}\overline{\bigcup_{s\ge
t}\varphi(s,\theta_{-s}\omega)D(\theta_{-s}\omega)}.
\]
\end{definition}

\begin{definition}
For given two random sets $D,A$, we say $A$ {\em (pull-back)
attracts} $D$ if
\[
\lim_{t\rightarrow\infty}d(\varphi(t,\theta_{-t}\omega)D(\theta_{-t}\omega)|A(\omega))=0
\]
holds almost surely, where $d(A|B)$ stands for the Hausdorff
semi-metric between two sets $A,B$, i.e. $d(A|B):={\rm sup}_{x\in
A}{\rm inf}_{y\in B}d_X(x,y)$; and we say $A$ {\em attracts $D$ in
probability} or {\em weakly attracts} $D$ if
\[
\mathbb P-\lim_{t\rightarrow\infty}d(\varphi(t,\omega)
D(\omega)|A(\theta_t\omega))=0.
\]
\end{definition}

By the measure preserving of $\theta_t$, it is clear that pull-back
attraction implies weak attraction.

\begin{remark}\label{omega}\rm
(i) Clearly $x\in\Omega_D(\omega)$ if and only if there exist
sequences $t_n\rightarrow\infty$ and $x_n\in D(\theta_{-t_n}\omega)$
such that $\varphi(t_n,\theta_{-t_n}\omega)x_n\rightarrow x$ as $n\rightarrow\infty$.\\
(ii) For any two given random sets $D_1$ and $D_2$, we have
\[
\Omega_{D_1\cup
D_2}(\omega)=\Omega_{D_1}(\omega)\cup\Omega_{D_2}(\omega)
\]
almost surely. By the definition of omega-limit sets, we clearly
have
\[
\Omega_{D_1\cup
D_2}(\omega)\supset\Omega_{D_1}(\omega)\cup\Omega_{D_2}(\omega),
\]
so we only need to show the converse inclusion holds. To see this,
for arbitrary $x\in \Omega_{D_1\cup D_2}(\omega)$, there exist
sequences $t_n\rightarrow\infty$ and  $x_n\in
D_1(\theta_{-t_n}\omega)\cup D_2(\theta_{-t_n}\omega)$ such that
$\varphi(t_n,\theta_{-t_n}\omega)x_n\rightarrow x$ as
$n\rightarrow\infty$. Hence there exists a subsequence such that
\[
x_{n_k}\in D_1(\theta_{-t_{n_k}}\omega)~~{\rm or}~~ x_{n_k}\in
D_2(\theta_{-t_{n_k}}\omega)
\]
holds for all $k=1,2,\ldots$ and
$\varphi(t_{n_k},\theta_{-t_{n_k}}\omega)x_{n_k}\rightarrow x$,
$k\rightarrow\infty$. That is, $x\in\Omega_{D_1}(\omega)$ or $x\in
\Omega_{D_2}(\omega)$. Therefore, $\Omega_{D_1\cup
D_2}(\omega)\subset\Omega_{D_1}(\omega)\cup\Omega_{D_2}(\omega)$.\\
(iii) If a random closed set $E$ pull-back attracts $D$, then
$\Omega_D\subset E$ almost surely. Indeed, if this is false, i.e.
the set $\hat\Omega:=\{\omega|~\Omega_D(\omega)\not\subset
E(\omega)\}$ has positive probability, for any
$\omega\in\hat\Omega$, assuming that $x\in\Omega_D(\omega)\backslash
E(\omega)$, then there exist sequences $t_n\rightarrow\infty$ and
$x_n\in D(\theta_{-t_n}\omega)$ such that
$\varphi(t_n,\theta_{-t_n}\omega)x_n\rightarrow x$ as
$n\rightarrow\infty$. Hence by the definition of Hausdorff
semi-metric and the fact that $E$ pull-back attracts $D$ we have
\begin{align*}
0<d(\{x\}|E(\omega))&=
\lim_{n\rightarrow\infty}d(\varphi(t_n,\theta_{-t_n}\omega)x_n|E(\omega))\\
&\le
\lim_{n\rightarrow\infty}d(\varphi(t_n,\theta_{-t_n}\omega)D(\theta_{-t_n}\omega)|E(\omega))\\
&=0,
\end{align*}
a contradiction.
\end{remark}

The following proposition comes from \cite{rud}, which gives a
relation between an $\bar{\mathcal F}^\nu$-measurable function and
an $\mathcal F$-measurable one.
\begin{proposition}\label{rud}
Assume that $\nu$ is a positive measure on the measurable space
$(X,\mathcal F)$. Denote $\bar{\mathcal F}^\nu$ the completion of
the $\sigma$-algebra $\mathcal F$ with respect to the measure $\nu$.
If $f$ is an $\bar{\mathcal F}^\nu$-measurable function, then there
exists an $\mathcal F$-measurable function $g$ such that $f=g$
$\nu$-a.e.
\end{proposition}

\section{Conley decomposition for random semiflow}

\subsection{Characterization of random chain recurrent set by random attractors (an observation)}

In this subsection, through simple observation, we obtain that, for
random semiflow, the characterization of random chain recurrent set
by random attractors holds similar to that in random flow case
presented in \cite{Liu,Liu2}.

\begin{definition}{\rm (\cite{Liu2})}
(i) Assume that $\epsilon>0$ is a random variable. A random open set
$U$ is called {\em $\epsilon$-absorbing} if there exists a random
variable $T>0$ such that $U$ contains the $\epsilon$-neighborhood of
$U_T(\omega):=\overline{\bigcup_{t\ge
T}\varphi(t,\theta_{-t}\omega)U(\theta_{-t}\omega)}$, i.e.
\[
B_{\epsilon}(U_T(\omega))\subset U(\omega).
\]
And we call a random open set $U$ {\em absorbing} if it is
$\epsilon$-absorbing for some random variable $\epsilon>0$.\\
(ii) An invariant random closed set $A$ is called an {\em (local)
attractor} if there exists an absorbing neighborhood $U$ of $A$ such
that $A(\omega)=\Omega_U(\omega)$. And we call
\[
B(A,U)(\omega):=\{x|~\varphi(t,\omega)x\in U(\theta_t\omega)~{\rm
for~ some}~t\ge 0\}
\]
the {\em basin of attraction of $A$ with respect to $U$}.
\end{definition}

\begin{definition}\label{def}
Assume that $A$ is an attractor with a random absorbing neighborhood
$U$ and the basin of attraction of $A$ with respect to $U$,
$B(A,U)$. Then we call
\[
R(\omega):=X\backslash B(A,U)(\omega)
\]
the {\em repeller corresponding to $A$ with respect to $U$}, and
call $(A,R)$ an {\em attractor-repeller pair of $\varphi$ (with
respect to $U$)}.
\end{definition}

\begin{remark}\label{rem}\rm
(i) Generally speaking, a random attractor is only forward
invariant, not necessarily invariant. In particular, when $\overline
U$ is a random compact set or $\varphi$ is a random ``flow" instead
of ``semiflow", the random attractor is invariant, see Lemma 3.2,
Proposition 3.6 and Remark 3.7 in \cite{Cra1} for details.\\
(ii) The repeller corresponding to $A$, $R$, is also forward
invariant. In fact, if there exists some $x_0\in R(\omega)$ and
$t_0>0$ such that $\varphi(t_0,\omega)x_0\notin
R(\theta_{t_0}\omega)$, i.e. $\varphi(t_0,\omega)x_0\in
B(A,U)(\theta_{t_0}\omega)$, then by the definition of $B(A,U)$ we
have
\[
\varphi(t_1,\theta_{t_0}\omega)\circ\varphi(t_0,\omega)x_0\in
U(\theta_{t_1}\circ\theta_{t_0}\omega)
\]
for some $t_1\ge 0$, that is
\[
\varphi(t_0+t_1,\omega)x_0\in U(\theta_{t_0+t_1}\omega).
\]
Hence we have $x_0\in B(A,U)(\omega)$, a contradiction.\\
(iii) It is easy to see that Lemma 5.2 in \cite{Liu2} also holds
when $\varphi$ is a random semiflow, i.e. for any given random
attractor $A$ there exists a forward invariant absorbing
neighborhood $U$ of $A$ such that $\Omega_U(\omega)=A(\omega)$. In
this case, it follows immediately that $B(A,U)$ is a forward
invariant random open set by the the definition of $B(A,U)$ and the
forward invariance of $U$.
\end{remark}

Similar to the proofs of Lemmas 3.4, 3.6 and 3.7 in \cite{Liu}, we
can obtain that Theorem 6 in \cite{Liu2} also holds for random
semiflow. In fact, the proofs of these lemmas are the same for both
random semiflow and random flow. Hence by (iii) of Remark \ref{rem}
we have the following

\begin{theorem}
Assume that $X$ is a Polish space and $\varphi$ is a random semiflow
on $X$. Assume that $U$ is a forward invariant random absorbing set,
$A$ is the random attractor determined by $U$ and $R$ is the random
repeller corresponding to $A$ with respect to $U$, then
\begin{equation}\label{mr1}
\mathcal{CR}_\varphi=\bigcap (A\cup R) ~~{\rm almost~surely},
\end{equation}
where $\mathcal{CR}_\varphi$ denotes the random chain recurrent set
of $\varphi$ (see \cite{Liu} for the definition for random flow;
this definition is also applicable to random semiflow), and the
intersection is taken over all forward invariant random absorbing
sets.
\end{theorem}

\begin{remark}\label{con}\rm
In above theorem, we characterize the random chain recurrent set by
random attractors for random semiflow on $X$. It is easy to see that
the corresponding result also holds when RDS $\varphi$ is restricted
to an invariant random compact set.
\end{remark}

\subsection{Complete Lyapunov function for random semiflow}

Firstly for arbitrary random attractor $A$, we will show that we can
construct a Lyapunov function $\hat l$ for $\varphi$ with respect to
$A$, see Lemmas \ref{l} and \ref{strict} for the construction.

\begin{lemma}\label{l}
Assume that $\mu$ is a positive measure on the product measurable
space $(\Omega\times X,\mathcal F\times\mathcal B(X))$, $A$ is a
random attractor with a forward invariant random absorbing
neighborhood $U$, and $R$ is the random repeller corresponding to
$A$ with respect to $U$. Then there exists an $\mathcal
F\times\mathcal B(X)$-measurable function $\tilde l(\omega,x)$ such
that
\[
\tilde l(\omega,x)=\left\{
\begin{array}{ll}
  0, & x\in A(\omega),\\
  1, & x\in R(\omega),\\
  0<\tilde l(\omega,x)\le 1, & x\in X\backslash(A(\omega)\cup
  R(\omega)).
\end{array}
\right.
\]
Moreover, $\tilde l(\omega,x)$ is non-increasing $\mu$-a.e. along
the orbits of the {\em skew-product semiflow $\Theta$ corresponding
to $\varphi$}, i.e.
$\Theta_t(\omega,x):=(\theta_t\omega,\varphi(t,\omega)x)$.
\end{lemma}

\noindent{\bf Proof.} Let
\[
\phi(\omega,x):=\frac{{\rm dist}_X(x,A(\omega))}{{\rm
dist}_X(x,A(\omega))+{\rm dist}_X(x,R(\omega))}.\] It is clear that
\[
\phi(\omega,x)=\left\{
\begin{array}{ll}
 0, &~  x\in A(\omega),\\
 1, &~ x\in R(\omega).
\end{array}
\right.
\]
Let
\[
l(\omega,x)=\sup_{t\geq 0}\phi(\theta_t\omega,\varphi(t,\omega)x).
\]
By the forward invariance of $A$ and $R$, we have
\begin{equation}\label{ly}
l(\omega,x)=\left\{
\begin{array}{ll}
  0, & x\in A(\omega),\\
  1, & x\in R(\omega),\\
  0<l(\omega,x)\le 1, & x\in X\backslash(A(\omega)\cup R(\omega)).
\end{array}
\right.
\end{equation}
By the definition of $l(\omega,x)$, it is clear that
\[
l(\theta_t\omega,\varphi(t,\omega)x)\le l(\omega,x),~\forall t>0,
\]
i.e. $l(\omega,x)$ is non-increasing along the orbits of $\Theta$.

Next we prove the measurability of $l(\omega,x)$. For $\forall
a\in\mathbb R^+$, we have
\begin{align*}
&\{ (\omega,x)|~l(\omega,x)>a\}\\
=&\{(\omega,x)|~\sup_{t\ge
0}\phi(\theta_t\omega,\varphi(t,\omega)x)>a\}\\
=& \Pi_{\Omega\times
X}\{(t,\omega,x)|~\phi(\theta_t\omega,\varphi(t,\omega)x)>a,~t\ge
0\},
\end{align*}
where $\Pi_{\Omega\times X}$ stands for the canonical projection of
$\mathbb R^+\times\Omega\times X$ to $\Omega\times X$. By the
measurability of the maps $\phi:(\omega,x)\mapsto\phi(\omega,x)$,
$\theta:(t,\omega)\mapsto\theta_t\omega$,
$\varphi:(t,\omega,x)\mapsto \varphi(t,\omega)x$, we know that the
map
\[
(t,\omega,x)\mapsto\phi(\theta_t\omega,\varphi(t,\omega)x)\] is
$\mathcal B(\mathbb R^+)\times\mathcal F\times\mathcal
B(X)$-measurable. Hence we obtain that
\[\{(t,\omega,x)|~\phi(\theta_t\omega,\varphi(t,\omega)x)>a,t\ge 0\}
\in\mathcal B(\mathbb R^+)\times\mathcal F\times\mathcal B(X).\] By
the projection theorem (see Proposition 2.3 in \cite{Liu}, which is
originated from \cite{Cas}) we have
\[\{ (\omega,x)|~l(\omega,x)>a\}\in(\mathcal F\times\mathcal B(X))^u,\]
that is, $l(\omega,x)$ is $(\mathcal F\times\mathcal
B(X))^u$-measurable, where the superscript $u$ denotes the universal
$\sigma$-algebra. By Proposition \ref{rud}, there exists an
$\mathcal F\times\mathcal B(X)$-measurable function $\tilde
l(\omega,x)=l(\omega,x)$ $\mu$-a.e. This completes the proof of the
lemma. \hfill $\Box$

\begin{remark}\label{remark1}\rm
The function $\tilde l(\omega,x)$ obtained in Lemma \ref{l} has two
shortcomings:
\begin{itemize}
    \item $0<\tilde l(\omega,x)\le 1$ when $x\in X\backslash(A(\omega)\cup R(\omega))$,
           not the desired $0<\tilde l(\omega,x)<1$;
    \item $\tilde l(\omega,x)$ is non-increasing along the orbits
          of $\Theta$ when $x\in X\backslash(A(\omega)\cup R(\omega))$, but
          not necessarily strictly decreasing.
\end{itemize}
\end{remark}

In fact, in the following lemma we will show that the two
shortcomings mentioned in Remark \ref{remark1} can be partially
improved, that is, $0<\tilde l(\omega,x)<1$ and strict decreasing
along the orbits of $\Theta$ when $x\in X\backslash(A(\omega)\cup
R(\omega))$ can hold in weaker sense---in probability.

\begin{lemma}\label{strict}
Assume that $A$ is a random attractor, and $R$ is the corresponding
random repeller of $A$ with respect to some random absorbing
neighborhood $U$. Then there exists a {\em Lyapunov function} $\hat
l(\omega,x)$ which has the same properties as $\tilde l(\omega,x)$
constructed in Lemma \ref{l}. Furthermore, for arbitrary
$\epsilon>0$ and arbitrary random variable  $x\in X\backslash(A\cup
R)$, there exists an $\Omega_\epsilon\subset\Omega$ with the measure
$\mathbb P(\Omega_\epsilon)<\epsilon$ such that on
$\Omega\backslash\Omega_\epsilon$, $\hat l(\omega,x)$ has the
following two fine properties:
\begin{eqnarray}
    & 0<\hat l(\omega,x(\omega))<1;\label{pro1}\\
& \hat l(\theta_t\omega,\varphi(t,\omega)x(\omega))<\hat
l(\omega,x(\omega)),
          ~\forall t>0.\label{pro2}
\end{eqnarray}
\end{lemma}

\noindent{\bf Proof.} Let
\begin{equation}\label{lya}
\hat l(\omega,x):=\frac12\left[\tilde
l(\omega,x)+\int_0^{+\infty}e^{-t}\tilde
l(\theta_{t}\omega,\varphi(t,\omega)x){\rm d}t\right].
\end{equation}
It is clear that $\hat l$ takes value $0$ on $A$ and takes value $1$
on $R$, respectively, by the forward invariance of $A$ and $R$. By
the monotonicity of $\tilde l$ along the orbits of $\varphi$, it
follows that $\hat l$ has the same monotonicity and that
\[
\int_0^{+\infty}e^{-t}\tilde
l(\theta_{t}\omega,\varphi(t,\omega)x){\rm
d}t\le\int_0^{+\infty}e^{-t}{\rm d}t\cdot\tilde l(\omega,x)=\tilde
l(\omega,x).
\]
Hence we have $\displaystyle\frac12\tilde l(\omega,x)\le\hat
l(\omega,x)\le\tilde l(\omega,x)$, which implies that $0<\hat
l(\omega,x)\le 1$ when $x\in X\backslash(A(\omega)\cup R(\omega))$.
That is, $\hat l(\omega,x)$ has the same properties as $\tilde
l(\omega,x)$.

For arbitrary random variable $x\in B(A,U)$, we have
\begin{equation}\label{limit}
\mathbb P-\lim_{t\rightarrow\infty}\tilde
l(\theta_t\omega,\varphi(t,\omega)x(\omega))=0.
\end{equation}
In fact, by the definition of $B(A,U)$, there exists a $t=t(\omega)$
such that $\varphi(t,\omega)x(\omega)\in U(\theta_{t}\omega)$. Then
for $\forall\epsilon>0$ we have
\begin{align}\label{8-1}
& \lim_{s\rightarrow\infty}\mathbb P\{\omega|~{\rm
dist}_X(\varphi(s,\omega)x(\omega),A(\theta_s\omega))>\epsilon\}\nonumber\\
=& \lim_{s\rightarrow\infty}\mathbb P\{\omega|~{\rm
dist}_X(\varphi(s+t,\omega)x(\omega),A(\theta_{s+t}\omega))>\epsilon\}\nonumber\\
=& \lim_{s\rightarrow\infty}\mathbb P\{\omega|~{\rm
dist}_X(\varphi(s,\theta_{t}\omega)\circ
\varphi(t,\omega)x(\omega),A(\theta_s\circ\theta_{t}\omega))>\epsilon\}\\
\le & \lim_{s\rightarrow\infty}\mathbb
P\{\omega|~d(\varphi(s,\theta_{t}\omega)U(\theta_{t}\omega)|
A(\theta_s\circ\theta_{t}\omega))>\epsilon\}\nonumber\\
=& 0,\nonumber
\end{align}
where the inequality holds by the definition of Hausdorff
semi-metric and the last equality holds by the fact $\Omega_U=A$ and
the measure preserving of $\theta_t$.

Now we will prove that (\ref{pro1}) and (\ref{pro2}) hold. The idea
of the proof is borrowed from \cite{Ar2}.

\indent If (\ref{pro1}) does not hold, i.e. there existing some
$\epsilon_0>0$ (without loss of generality, suppose $\epsilon_0\le
1$) and some random variable $x_0\in B(A,U)\backslash A$, for
arbitrary $\Omega_{\epsilon_0}$ with $\mathbb
P(\Omega_{\epsilon_0})< \epsilon_0$, there exists an
$\tilde\Omega\subset\Omega\backslash\Omega_{\epsilon_0}$ with
$\mathbb P(\tilde\Omega)>0$ such that on $\tilde\Omega$
\[
\hat l(\omega,x_0(\omega))=1.
\]
Then by (\ref{lya}), we have
\[\tilde l(\theta_t\omega,\varphi(t,\omega)x_0(\omega))=1,
~\forall\omega\in\tilde\Omega, ~{\rm Leb-almost~all}~t\in\mathbb
R^+,\] a contradiction to (\ref{limit}).

\indent If (\ref{pro2}) does not hold, similar to the proof of
(\ref{pro1}), there existing some $\epsilon_0>0$ and some random
variable $x_0\in B(A)\backslash A$, for arbitrary
$\Omega_{\epsilon_0}$ with $\mathbb P(\Omega_{\epsilon_0})<
\epsilon_0$, there exists an
$\tilde\Omega\subset\Omega\backslash\Omega_{\epsilon_0}$ with
$\mathbb P(\tilde\Omega)>0$ such that for arbitrary
$\omega\in\tilde\Omega$, there exists $t_0=t_0(\omega)>0$ satisfying
\[
\hat l(\theta_{t_0}\omega,\varphi(t_0,\omega)x_0(\omega))=\hat
l(\omega,x_0(\omega)).
\]
By the monotonicity of $\tilde l$ along the orbits of $\Theta$, we
have
\begin{equation}\label{A1}
\tilde l(\theta_s\omega,\varphi(s,\omega)x_0(\omega))=\tilde
l(\omega,x_0(\omega))>0 ~{\rm for~all}~0\le s\le t_0,
\end{equation}
and
\begin{equation*}
\tilde
l(\theta_{s+t_0}\omega,\varphi(s+t_0,\omega)x_0(\omega))=\tilde
l(\theta_{s}\omega,\varphi(s,\omega)x_0(\omega)) ~{\rm
for~Leb-almost ~all}~s\ge 0.
\end{equation*}
Hence
\begin{equation}\label{A2}
\tilde
l(\theta_{nt_0+s}\omega,\varphi(nt_0+s,\omega)x_0(\omega))=\tilde
l(\theta_s\omega,\varphi(s,\omega)x_0(\omega))
\end{equation}
for all $n\in\mathbb N$ and for Leb-almost all $s\ge 0$. Therefore
for each $\omega\in\tilde\Omega$, $\exists \tau=\tau(\omega)\ge 0$
such that for which both (\ref{A1}) and (\ref{A2}) hold, i.e. we
have
\begin{equation}\label{A3}
\tilde
l(\theta_{nt_0+\tau}\omega,\varphi(nt_0+\tau,\omega)x_0(\omega))=\tilde
l(\omega,x_0(\omega))>0,~~\forall n\in\mathbb N.
\end{equation}
Letting $n\rightarrow\infty$ in (\ref{A3}), we obtain a
contradiction to (\ref{limit}). This terminates the proof of the
lemma. \hfill $\Box$

By completely similar to \cite{Liu2} in random flow case, we can
also construct complete Lyapunov function for random semiflow;
furthermore, we can discuss chain transitive components completely
similar to random flow case. In fact, it is clear that semiflow or
flow is not relevant in these steps. Hence here we only state
associated results and omit details of the proof.

\begin{theorem}\label{LY}
Assume that $X$ is a Polish space and $\varphi$ is a random semiflow
on $X$. Then there exists a complete Lyapunov function
$L:\Omega\times X\rightarrow\mathbb R^+$ for $\varphi$ with the
following
properties:\\
(i) $L$ is an $\mathcal F\times\mathcal B(X)$-measurable function; \\
(ii) $L(\theta_t\omega,\varphi(t,\omega)x)\le L(\omega,x)$
               for $\forall t>0$ $\mu$-a.e., recalling that $\mu$ is defined in Lemma \ref{l};\\
(iii) $L(\theta_t\omega,\varphi(t,\omega)x)=L(\omega,x)$ for
               $\forall t>0$ when $x\in\mathcal{CR}_\varphi(\omega)$;\\
(iv) If the random variable $x$ is completely random non-chain
recurrent, i.e. $x(\omega)\in
X\backslash\mathcal{CR}_\varphi(\omega)$ $\mathbb P$-a.s., then for
arbitrary $\epsilon>0$ there exists an
          $\Omega_\epsilon\subset\Omega$ satisfying $\mathbb P(\Omega_\epsilon)<\epsilon$
          such that for arbitrary
          $\omega\in\Omega\backslash\Omega_\epsilon$, the
          following holds:
          \[L(\theta_t\omega,\varphi(t,\omega)x(\omega))<L(\omega,x(\omega)),~\forall
          t>0.\]
(v) The range of $L$ on $\mathcal{CR}_\varphi(\omega)$ is a compact
nowhere dense
               subset of $[0,1]$;\\
(vi) $L$ separates different random chain transitive components of
$\varphi$;\\
(vii) If $C$ and $C'$ are distinct random chain transitive
          components of $\varphi$ with the property that for
          arbitrary random variables
          $\epsilon,T>0$
          there is an $\epsilon$-$T$-chain from
          $C$ to $C'$ $\mathbb P$-a.s. then $L(\Omega,C)>L(\Omega,C')$.
\end{theorem}

\begin{remark}\rm
Among the properties of complete Lyapunov function for random
semiflow, (ii) and (iv) are weaker than that of random flow case
obtained in \cite{Liu2}. The reason is that the Lyapunov function
for attractor-repeller pair constructed in Lemma \ref{strict} has
weaker properties.
\end{remark}

\section{Morse decomposition of global random attractors}

Global random attractors were introduced by Crauel and Flandoli
\cite{Cra1}, Schmalfuss \cite{Sch}, and were studied for many SDEs,
see \cite{CDF,FS,Sche,Sch2}, among others. First let us recall the
definition of global random attractor.

\begin{definition}(\cite{Cra1})\label{attractor}
Assume that $\varphi$ is a random semiflow on a Polish space $X$,
then a random compact set $A$ is called a  {\em global random
attractor} for $\varphi$ if
\begin{itemize}
  \item $A$ is invariant, i.e.
\begin{equation}\label{inva}
\varphi(t,\omega)A(\omega)=A(\theta_t\omega),~\forall t\ge 0
\end{equation}
for almost all $\omega\in\Omega$;

  \item $A$ pull-back attracts every bounded deterministic
set, i.e. for any bounded deterministic set $B\subset X$, we have
\begin{equation}\label{att}
\lim_{t\rightarrow\infty}d(\varphi(t,\theta_{-t}\omega)B|A(\omega))=0
\end{equation}
almost surely.
\end{itemize}
\end{definition}

The global random attractor for RDS $\varphi$ is the {\em minimal}
random compact set which attracts all the bounded deterministic sets
and it is the {\em largest} random compact set which is invariant in
the sense of (\ref{inva}), see \cite{CDF} for details. The global
random attractor is uniquely determined by attracting {\em
deterministic compact} sets, see \cite{Cr} for details.

We now introduce ``backward orbit" for random semiflow:
\begin{itemize}
  \item For fixed $\omega$ and $x$, a mapping $\sigma_\cdot(\omega):\mathbb R^-\rightarrow X$
         is called a {\em backward orbit of $\varphi$ through $x$ driven by $\omega$} if it
         satisfies the cocycle property:
         \[
         \sigma_0(\omega)=x,~
         \sigma_{t+s}(\omega)=\varphi(s,\theta_t\omega)\circ\sigma_t(\omega)~{\rm for}~
         \forall t\le 0, s\ge 0, t+s\le0.
         \]
  \item Let $\mathcal M$ denote the set of all $X$-valued random
  variables and $x\in\mathcal M$. A mapping $\sigma:\mathbb R^-\rightarrow\mathcal M$ is
  called a {\em backward orbit of $\varphi$ through $x$} if for all
  $\omega\in\Omega$, the following cocycle property holds:
  \[
\sigma_0(\omega)=x(\omega),~
         \sigma_{t+s}(\omega)=\varphi(s,\theta_t\omega)\circ\sigma_t(\omega)~{\rm for}~
         \forall t\le 0, s\ge 0, t+s\le0.
  \]
\end{itemize}
Also we can introduce ``entire orbit" for random semiflow:
\begin{itemize}
  \item For fixed $\omega$ and $x$, a mapping $\sigma_\cdot(\omega):\mathbb R\rightarrow X$
         is called an {\em entire orbit of $\varphi$ through $x$ driven by $\omega$} if it
         satisfies the cocycle property:
         \[
         \sigma_0(\omega)=x,~
         \sigma_{t+s}(\omega)=\varphi(s,\theta_t\omega)\circ\sigma_t(\omega)~{\rm for}~
         \forall t\in\mathbb R, s\ge 0.
         \]
  \item Let $x\in\mathcal M$. A mapping $\sigma:\mathbb R\rightarrow\mathcal M$ is
  called an {\em entire orbit of $\varphi$ through $x$} if for all
  $\omega\in\Omega$, the following cocycle property holds:
  \[
\sigma_0(\omega)=x(\omega),~
         \sigma_{t+s}(\omega)=\varphi(s,\theta_t\omega)\circ\sigma_t(\omega)~{\rm for}~
         \forall t\in\mathbb R, s\ge 0.
  \]
\end{itemize}
By the cocycle property of $\sigma$, it is clear that for arbitrary
$t\ge 0$ and $\omega\in\Omega$ we have
\[
\sigma_t(\omega)=\varphi(t,\omega)\circ\sigma_0(\omega).
\]
That is, when an entire orbit $\sigma$ of $\varphi$ is restricted to
$\mathbb R^+$ (called {\em forward orbit}), then it coincides with
the orbit of $\varphi$, which is the same as the deterministic case.

We can give an alternative definition of (forward, backward)
invariant sets for random semiflow:
\begin{itemize}
\item A random set $D$ is called {\em forward invariant} if
  $D=D_\varphi^+$ almost surely, where
  \[
  D_\varphi^+(\omega):=\{x|\varphi(t,\omega)x\in D(\theta_t\omega)~{\rm
  for~all}~t\ge0\};
  \]
\item A random set $D$ is called {\em backward invariant} if
  $D=D_\varphi^-$ almost surely, where
  \[
   D_\varphi^-(\omega):=\{x|{\rm \exists ~a~ backward~ orbit~\sigma~ in}~D~{\rm
  through}~x,~{\rm i.e.}~
  \sigma_t(\omega)\in D(\theta_t\omega),\forall t\le 0\};
  \]
\item A random set $D$ is called {\em invariant} if
  $D=D_\varphi$ almost surely, where
  \[
 D_\varphi(\omega):=\{x|{\rm \exists ~an~ entire~
  orbit~\sigma~in~}D~{\rm
  through}~x,~{\rm i.e.}~\sigma_t(\omega)\in D(\theta_t\omega),\forall t\in\mathbb R\}.
  \]
\end{itemize}
\begin{remark}\label{rem4.1}\rm
(i) Clearly a random set $D$ is invariant if and only if it is both
forward invariant and backward invariant. It is easy to verify that
forward invariant and invariant sets defined above coincide with
that of Definition \ref{invariant}, which is convenient for us to
choose appropriate definition in the sequel.\\
(ii) If $D_1$ and $D_2$ are forward invariant, then clearly $D_1\cup
D_2$ and $D_1\cap D_2$ are forward invariant; If $D_1$ and $D_2$ are
invariant, then clearly $D_1\cup D_2$ is invariant, while $D_1\cap
D_2$ is not necessarily invariant (since $D_1\cap D_2$ is not
necessarily backward invariant), which differs from random flow
case, see page 35 in \cite{Ar1}.
\end{remark}

\begin{lemma}\label{exten}
Assume that $D$ is a forward invariant random compact set, then for
any point on $\Omega_D$ there exists a backward orbit lying on
$\Omega_D$ through this point.
\end{lemma}

\noindent{\bf Proof.} By the definition of omega-limit set we know
that, for any given $x\in\Omega_D(\omega)$, there exist sequences
$t_n\rightarrow+\infty$, $x_n\in D(\theta_{-t_n}\omega)$ such that
$\varphi(t_n,\theta_{-t_n}\omega)x_n\rightarrow x$,
$n\rightarrow\infty$. For arbitrary $k\in\mathbb Z^-$, there exists
an $N_0$ such that the set
$\{\varphi(t_n+k,\theta_{-t_n}\omega)x_n|n\ge N_0\}$ is pre-compact,
i.e. the closure of $\{\varphi(t_n+k,\theta_{-t_n}\omega)x_n|n\ge
N_0\}$ is compact. In fact, by the forward invariance of $D$, when
$t_n+k\ge 0$, we have
\[
\{\varphi(t_n+k,\theta_{-t_n}\omega)x_n|n\ge N_0\}\subset
D(\theta_k\omega).
\]
By taking a subsequence,
$\lim_{n\rightarrow\infty}\varphi(t_n+k,\theta_{-t_n}\omega)x_n$
exists and denote by $\tilde x_k$ the limit. Since
$D(\theta_k\omega)$ is compact, we have $\tilde x_k\in
D(\theta_k\omega)$. For $k\le t\le k+1$, $k\in\mathbb Z^-$, let
$\sigma_t(\omega)=\varphi(t-k,\theta_{k}\omega)\tilde x_k$. Clearly
$\sigma$ is a backward orbit of $\varphi$ through $x$. Moreover the
backward orbit obtained in this way lies on $\Omega_D$. In fact, by
the definition of $\Omega_D$, it is clear that $\tilde
x_k\in\Omega_D(\theta_k\omega)$. Hence by the invariance of
$\Omega_D$, we have $\sigma_t(\omega)\in\Omega_D(\theta_t\omega)$,
$t\le 0$. \hfill $\Box$

\begin{corollary}\label{cor}
Assume that $D$ is an invariant random compact set, then for any
point on $D$, there exists a backward orbit lying on $D$ through
this point.
\end{corollary}
\noindent{\bf Proof.} It follows immediately from Lemma \ref{exten}.
\hfill $\Box$

In contrast to Lemma \ref{exten}, given a forward invariant random
compact set $D$, a natural question is, for any random variable
$x\in\Omega_D$, does there exist a backward orbit lying on
$\Omega_D$ through $x$? The answer is yes, see the following lemma.

\begin{lemma}
Assume that $D$ is a forward invariant random compact set, then for
any random variable on $\Omega_D$ there exists a backward orbit
lying on $\Omega_D$ through this random variable.
\end{lemma}

\noindent{\bf Proof.} By Lemma \ref{exten} we know that, for given
$k\in\mathbb Z^-$ and for each $\omega$, there exists an $\tilde
x_k(\omega)\in\Omega_D(\theta_k\omega)$ from which we obtain a
backward orbit from time $k$ to time $0$ (present time). Hence we
need only to show that we can select appropriate $\tilde x_k$ such
that the map $\omega\mapsto\tilde x_k(\omega)$ is measurable. In
other words, we need to show $\sigma_k\in\mathcal M$, which implies
$\sigma_s\in\mathcal M$, $\forall k\le s\le0$. For arbitrary $t>0$,
denote $\varphi^{-1}(t,\omega)x$ the preimage of $x$ under
$\varphi$. Consider the skew-product semiflow $\Theta$ corresponding
to $\varphi$ (see Lemma \ref{l}) which is an $\mathcal
F\times\mathcal B(X)$-measurable mapping from $\Omega\times X$ to
itself for fixed $t\ge 0$. The preimage of $(\omega,x)$ under
$\Theta_t$ is
$\Theta_t^{-1}(\omega,x):=(\theta_{-t}\omega,\varphi^{-1}(t,\omega)x)$.
Since $\Omega_D$ is a random compact set, we have
\[
{\rm graph}(\Omega_D):=\{(\omega,x)|x\in\Omega_D(x)\}\in\mathcal
F\times\mathcal B(X),
\]
see page 59 of \cite{Cas} or Proposition 2.4 in \cite{Cra02}. Hence
we have $\Theta_t^{-1}({\rm graph}(\Omega_D))\in\mathcal
F\times\mathcal B(X)$ by the measurability of $\Theta_t$. It is
cleat that
\[
{\rm graph}(\Omega_D^t)=\Theta_t^{-1}({\rm graph}(\Omega_D)),
\]
where
\[
\Omega_D^t(\omega):=\varphi^{-1}(t,\omega)\Omega_D(\theta_t\omega).
\]
Therefore $\Omega_D^t$ is an $\mathcal F^u$-measurable random
compact set, see page 59 of \cite{Cas} or Proposition 2.4 in
\cite{Cra02} again. By Lemma 2.7 in \cite{Cra02} we know that we may
assume that $\Omega_D^t$ is an $\mathcal F$-measurable random
compact set. In particular, for given $k\in\mathbb Z^-$,
$\Omega_D^{-k}\cap\Omega_D$ is a nonempty random compact set by (v)
of Proposition 2.1 in \cite{Liu}. By the measurable selection
theorem (see Proposition 2.2 in \cite{Liu}), we can choose a random
variable $\tilde x_k\in\Omega_D^{-k}\cap\Omega_D$. This completes
the proof. \hfill $\Box$

So we have
\begin{corollary}\label{corol}
Assume that $D$ is an invariant random compact set, then for any
random variable on $D$, there exists a backward orbit lying on $D$
through this random variable.
\end{corollary}

Throughout this section, we use $S$ to denote the invariant random
compact set we will decompose, say, $S$ is a global random
attractor. By Corollaries \ref{cor} and \ref{corol}, for any point
(random variable) on $S$, there exists backward orbit lying on $S$
through this point (random variable). Afterwards, when we say
backward orbits, we refer those lying on $S$ unless otherwise stated
(since there may be backward orbit not lying on $S$ but lying on the
entire state space --- $X$).

\begin{definition}
An invariant random compact set $A\subset S$ is called a {\em
(local) attractor} if there exists a random open neighborhood $U$ of
$A$ relative to $S$ such that $\Omega_U(\omega)=A(\omega)$. (By
Lemma 3.1 in \cite{Liu2}, without loss of generality, we can assume
that $U$ is forward invariant.) The basin of attraction of $A$ is
defined by
\[
B(A)(\omega):=\{x\in S(\omega)|\varphi(t,\omega)x\in
U(\theta_t\omega)~{\rm for~some~} t\ge0\}
\]
and the {\em repeller $R$ corresponding to $A$} is defined by
\[
R(\omega)=S(\omega)\backslash B(A)(\omega).
\]
$(A, R)$ is called an {\em attractor-repeller pair}.
\end{definition}

Note that since $S$ is a random compact set, by Lemma 3.2 in
\cite{Liu} (the proof of Lemma 3.2 is also applicable here), $B(A)$
is independent of the choice of $U$.

\begin{lemma}\label{lem2}
Assume that $(A,R)$ is an attractor-repeller pair in $S$, then $A$,
$B(A)$, and $R$ are invariant random sets.
\end{lemma}

\noindent{\bf Proof.} (i) The invariance of $A$ follows immediately from its definition. \\
(ii) The forward invariance of $B(A)$ follows directly from the
definition of $B(A)$ and the forward invariance of $U$. For
arbitrary $x\in B(A)(\omega)$, if for any backward orbit $\sigma$
through $x$, there exists some $t_0<0$ such that $\sigma_{t_0}\in
R(\theta_{t_0}\omega)$. By the definition of $R$, we know that any
point in $R$ can not enter into $B(A)$ in positive time, so we
obtain that $R$ is forward invariant. Therefore,
\[
\varphi(-t_0,\theta_{t_0}\omega)\sigma_{t_0}(\omega)=\sigma_0(\omega)=x\in
R(\omega),
\]
a contradiction. That is, $B(A)$ is backward invariant. \\
(iii) By (ii) we only need to show the backward invariance of $R$.
For arbitrary $x\in R(\omega)$, if for any backward orbit $\sigma$
through $x$, there exists some $t_0<0$ such that $\sigma_{t_0}\in
B(A)(\theta_{t_0}\omega)$. Then by the forward invariance of $B(A)$,
we have
\[
\varphi(-t_0,\theta_{t_0}\omega)\sigma_{t_0}(\omega)=\sigma_0(\omega)=x\in
B(A)(\omega),
\]
a contradiction. Hence $R$ is backward invariant. \hfill $\Box$

\begin{definition}\label{def-lim}
Assume that $x$ is a random variable in $S$, and $\sigma$ is an
entire orbit through $x$. Then the {\em omega-limit set
$\Omega_\sigma$} and the {\em alpha-limit set $\Omega_\sigma^*$} of
$\sigma$ are defined to be
\[
\Omega_\sigma(\omega):=\bigcap_{T\ge0}\overline{\bigcup_{t\ge
T}\{\sigma_t(\theta_{-t}\omega)\}}
\]
and
\[
\Omega_\sigma^*(\omega):=\bigcap_{T\ge0}\overline{\bigcup_{t\ge
T}\{\sigma_{-t}(\theta_{t}\omega)\}},
\]
respectively.
\end{definition}
It is clear that
\[
\Omega_\sigma(\omega):=\bigcap_{T\ge0}\overline{\bigcup_{t\ge
T}\{\varphi(t,\theta_{-t}\omega)x(\theta_{-t}\omega)\}},
\]
i.e. the omega-limit set of $\sigma$ only depends on the random
variable $x$, so $\Omega_\sigma$ can also be denoted by $\Omega_x$;
while the alpha-limit set depends on the entire orbit $\sigma$.
Clearly a point $y\in\Omega_\sigma(\omega)$ (respectively
$y\in\Omega_\sigma^*(\omega)$) if and only if there exist sequences
$t_n\rightarrow +\infty$ (respectively $t_n\rightarrow -\infty$) and
$y_n=\sigma_{t_n}(\theta_{-t_n}\omega)$ such that $y_n\rightarrow y$
as $n\rightarrow+\infty$.

\begin{lemma}\label{lem3}
Assume that $x$ is a random variable in $S$, and $\sigma$ is an
entire orbit through $x$. Then $\Omega_\sigma$ and $\Omega_\sigma^*$
are invariant random compact sets.
\end{lemma}

\noindent{\bf Proof.} The random variable $x$ can be regarded as a
random set consisting of just a single point, so $\Omega_\sigma$ is
an invariant random compact set.

For arbitrary $y\in\Omega_\sigma^*(\omega)$, there exist sequences
$t_n\rightarrow +\infty$ and $y_n=\sigma_{-t_n}(\theta_{t_n}\omega)$
such that $y_n\rightarrow y$ as $n\rightarrow+\infty$. For $s>0$, we
have
\begin{align*}
\varphi(s,\omega)y
&=\lim_{n\rightarrow+\infty}\varphi(s,\omega)\circ\sigma_{-t_n}(\theta_{t_n}\omega)\\
&=\lim_{n\rightarrow+\infty}\varphi(s,\omega)\circ\sigma_{-t_n}(\theta_{t_n-s}\circ\theta_s\omega)\\
&=\lim_{n\rightarrow+\infty}\sigma_{s-t_n}(\theta_{t_n-s}\circ\theta_s\omega)\\
&=\lim_{n\rightarrow+\infty}\sigma_{-\tau_n}(\theta_{\tau_n}\circ\theta_s\omega)~~({\rm
let~}t_n-s=\tau_n)\\
&\in\Omega_\sigma^*(\theta_s\omega),
\end{align*}
where the 1st equality holds by the continuity property of $\varphi$
with respect to $x$, the 3rd equality holds by the cocycle property
of $\sigma$, and the last inclusion relation holds by the definition
of $\Omega_\sigma^*$. This verifies the forward invariance of
$\Omega_\sigma^*$.

For arbitrary $y\in\Omega_\sigma^*(\theta_s\omega)$ with $s>0$,
there exist sequences $t_n\rightarrow +\infty$ and
$y_n=\sigma_{-t_n}(\theta_{t_n}\circ\theta_s\omega)$ such that
$y_n\rightarrow y$ as $n\rightarrow+\infty$. Then we have
\begin{align*}
y
&=\lim_{n\rightarrow+\infty}\sigma_{-t_n}(\theta_{t_n}\circ\theta_s\omega)\\
&=\lim_{n\rightarrow+\infty}\sigma_{-(\tau_n-s)}(\theta_{\tau_n-s}\circ\theta_s\omega)~~({\rm
let~}t_n+s=\tau_n)\\
&=\lim_{n\rightarrow+\infty}\varphi(s,\omega)\sigma_{-\tau_n}(\theta_{\tau_n-s}\circ\theta_s\omega)\\
&=\lim_{n\rightarrow+\infty}\varphi(s,\omega)\sigma_{-\tau_n}(\theta_{\tau_n}\omega)\\
&=\varphi(s,\omega)\lim_{n\rightarrow+\infty}\sigma_{-\tau_n}(\theta_{\tau_n}\omega)\\
&=\varphi(s,\omega)x,
\end{align*}
where the last two equalities hold by the pre-compactness of
$\{\sigma_{-\tau_n}(\theta_{\tau_n}\omega)|n\in\mathbb N\}$, and by
taking a subsequence we assume that the subsequence converges to
$x\in\Omega_\sigma^*$. This verifies
$\Omega_\sigma^*(\theta_s\omega)\subset\varphi(s,\omega)\Omega_\sigma^*(\omega)$.

Therefore, we have showed that
$\varphi(s,\omega)\Omega_\sigma^*(\omega)=\Omega_\sigma^*(\theta_s\omega)$,
hence completed the proof. \hfill $\Box$

\begin{lemma}\label{prop}
Assume that $x$ is a random variable with $\sigma$ being an entire
orbit
through $x$, and $A$ is a random attractor with $R$ being the corresponding repeller. Then we have:\\
(i) if $x\in R$ almost surely, then $\Omega_\sigma\subset R$ and $\Omega_\sigma^*\subset R$ almost surely;\\
(ii) if $x\in B(A)\backslash A$ almost surely, then
$\Omega_\sigma\subset A$ and $\Omega_\sigma^*\subset R$ almost
surely;\\
(iii) if $x\in A$ almost surely, then $\Omega_\sigma\subset A$
almost surely; if
    $\Omega_\sigma^*\subset A$ almost surely, then $\sigma$ lies on $A$ almost surely,
    i.e. for arbitrary $t\in\mathbb R$, we have $\sigma_t\subset A$ almost surely;  \\
(iv) if $x\in B(A)$ almost surely, then $\Omega_\sigma\subset A$
almost surely; if $x\in B(R):=S\backslash A$ almost surely, then
$\Omega_\sigma^*\subset R$ almost surely.
\end{lemma}

\noindent{\bf Proof.} (i) By the forward invariance of $R$, the
former is obvious. By the forward invariance of $B(A)$ we obtain
that all backward orbits through $x$ must lie on $R$, so by the
definition of $\Omega_\sigma^*$  we have $\Omega_\sigma^*\subset R$
almost surely.

(ii) The former follows directly from Lemma 4.3 in \cite{Liu3} (it
is clear that Lemma 4.3 also holds for random semiflow). Assume that
$U$ is a forward invariant random open neighborhood of $A$ relative
to $S$ such that $\Omega_U=A$ and let $V=S\backslash U$. For
arbitrary random variable $y\in V$, let $\sigma^y$ be a backward
orbit through $y$. Then by the forward invariance of $U$, $\sigma^y$
lies on $V$. Hence we have $\Omega_{\sigma^y}^*\subset V$ almost
surely. If $\Omega_{\sigma^y}^*\not\subset R$ with positive
probability, letting $R_1:=R\cup\Omega_{\sigma^y}^*$, then $R_1$ is
an invariant random compact set by Lemma \ref{lem3} and (ii) of
Remark \ref{rem4.1}. Then we can choose a random variable $z$ such
that $z\in R_1$ almost surely and $z\in R_1\backslash R$ with
positive probability. On one hand, $\Omega_{z}\subset R_1$ almost
surely by the invariance of $R_1$, which implies
\[
\mathbb
P-\lim_{t\rightarrow\infty}d(\varphi(t,\omega)z(\omega)|R_1(\theta_t\omega))=0.
\]
On the other hand $z\in B(A)$ with positive probability, which
implies that
\[
\lim_{t\rightarrow\infty}d(\varphi(t,\cdot)z(\omega)|A(\theta_t\cdot))=0
\]
with positive probability, a contradiction to the fact that $R_1\cap
A=\emptyset$ almost surely. Therefore, for arbitrary random variable
$y\in V$, we have $\Omega_{\sigma^y}^*\subset R$ almost surely.

Let $U_n:=\varphi(n,\theta_{-n}\omega)U(\theta_{-n}\omega)$,
$n\in\mathbb N$, then we have $U_{n+1}\subset U_n$ by the forward
invariance of $U$. Moreover, each $U_n$ is a forward invariant
random open neighborhood of $A$ relative to $S$ and
$\Omega_{U_n}=A$. Clearly we have
\[
A(\omega)=\lim_{n\rightarrow\infty}U_n(\omega).
\]
Letting $V_n=S\backslash U_n$, for arbitrary random variable $x\in
V_n$, we have $\Omega_{\sigma^x}^*\subset R$ almost surely by the
above argument. Since $n$ is arbitrary, for arbitrary random
variable in $S\backslash A$ with a backward orbit $\sigma$, we have
$\Omega_\sigma^*\subset R$ almost surely. This completes the proof
of (ii).

(iii) The former is trivial. Since $A$ is forward invariant, any
backward orbit through a random variable in $S\backslash A$ must lie
on $S\backslash A$. If there exists some $t_0\in\mathbb R$ such that
$\mathbb P\{\omega|\sigma_{t_0}(\omega)\not\subset
A(\theta_{t_0}\omega)\}=\delta>0$, then for all $s\le t_0$ we have
$\mathbb P\{\omega|\sigma_{s}(\omega)\not\subset
A(\theta_{s}\omega)\}\ge\delta$, i.e. $\mathbb
P\{\omega|\sigma_{s}(\omega)\subset S(\theta_{s}\omega)\backslash
A(\theta_{s}\omega)\}\ge\delta$. Then by the proof of (ii) it
follows that $\Omega_\sigma^*\subset R$ with positive probability, a
contradiction to the fact that $\Omega_\sigma^*\subset A$ almost
surely.

(iv) The former follows directly from (ii) and (iii), while the
later has been proved during the proof of (ii). \hfill $\Box$

\begin{definition} Assume that $(A_i,R_i)$ are
attractor-repeller pairs of $\varphi$ on the invariant random
compact set $S$ with
\[
\emptyset=A_0\varsubsetneq A_1\varsubsetneq\cdots\varsubsetneq
A_n=S~{\rm and}~S=R_0\varsupsetneq
R_1\varsupsetneq\cdots\varsupsetneq R_n=\emptyset.
\]
Then the family $D=\{M_i\}_{i=1}^{n}$ of invariant random compact
sets, defined by
\[
M_i=A_i\cap R_{i-1},~1\le i\le n,
\]
is called a {\em Morse decomposition of $S$}, and each $M_i$ is
called {\em Morse set}. If $D$ is a Morse decomposition, $M_D$ is
defined to be $\bigcup_{i=1}^nM_i$.
\end{definition}

\begin{remark}\rm
(i) For $i\neq j$, say $i<j$, $M_i\cap M_j=A_i\cap R_{i-1}\cap
A_j\cap R_{j-1}\subset A_i\cap R_{j-1}\subset A_i\cap
R_i=\emptyset$.\\
(ii) Each Morse set $M_i$ is invariant, which is trivial if
$\varphi$ is a random flow, but requires explanation in the case of
random semiflow. Clearly each $M_i$ is forward invariant. For any
point $x\in M_i(\omega)=A_i(\omega)\cap R_{i-1}(\omega)$, there
exists a backward orbit $\sigma$ in $A_i$ by the backward invariance
of $A_i$. By the forward invariance of $R^c_{i-1}=B(A_{i-1})$, any
backward orbit through a point in $R_{i-1}$ must lie on $R_{i-1}$.
Hence we have $\sigma$ lying on $M_i$, i.e. $M_i$ is backward
invariant.
\end{remark}

\begin{lemma}\label{lem1}
Assume that $A_1, A_2\subset S$ are two random attractors with
basins of attraction $B(A_1)$, $B(A_2)$, respectively. Assume that
$D$ is a random compact set satisfying $ D\subset B(A_1)\cup B(A_2)$
almost surely. Then $A_1\cup A_2$ pull-back attracts $D$.
\end{lemma}

\noindent{\bf Proof.} Denote
\[
\Omega_1=\{\omega|D(\omega)\subset
B(A_1)(\omega)\},~\Omega_2=\{\omega|D(\omega)\subset
B(A_2)(\omega)\backslash B(A_1)(\omega)\}.
\]
Choose a random compact set $D_1\subset B(A_1)$ almost surely
satisfying $D_1(\omega)=D(\omega)$ for $\omega\in\Omega_1$ and
choose $D_2\subset B(A_2)$ almost surely satisfying
$D_2(\omega)=D(\omega)$ for $\omega\in\Omega_2$. Then we have
$D\subset D_1\cup D_2$ almost surely. Therefore, by (ii) and (iii)
of Remark \ref{omega} and Lemma 4.3 in \cite{Liu3}, we obtain for
$\mathbb P$-almost all $\omega$
\[
\Omega_D(\omega)\subset\Omega_{D_1\cup
D_2}(\omega)=\Omega_{D_1}(\omega)\cup\Omega_{D_2}(\omega)\subset
A_1(\omega)\cup A_2(\omega).
\]
By the definition of omega-limit set, it is clear that $\Omega_D$
pull-back attracts $D$, so $A_1\cup A_2$ pull-back attracts $D$.
This completes the proof.  \hfill$\Box$

\begin{remark}\rm
It is obvious that the result of Lemma \ref{lem1} holds for finite
case, i.e. if the random compact set
$D\subset\bigcup_{i=1}^{n}B(A_i)$ almost surely, then
$\bigcup_{i=1}^{n}A_i$ pull-back attracts $D$.
\end{remark}

\begin{theorem}\label{Mor}
Assume that $D=\{M_i\}_{i=1}^{n}$ is a Morse decomposition of $S$,
determined by attractor-repeller pairs $(A_i,R_i)$, $i=1,\ldots,n$.
Then $M_D$ determines the limiting behavior of $\varphi$ on $S$.
Moreover, there are no ``cycles" between the Morse sets. More
precisely, we have:
\begin{description}
\item[(i)] For any random variable $x$ in $S$, there exists an entire orbit $\sigma$ through $x$
    such that $\Omega_\sigma\subset M_D$
    and $\Omega_\sigma^*\subset M_D$ almost surely.
\item[(ii)] If $\sigma$ is an entire orbit through the random variable
    $x$ satisfying that
    $\Omega_\sigma\subset M_p$ almost surely and $\Omega_\sigma^*\subset M_q$
    almost surely for some $1\le p,q\le
    n$, then $p\le q$;
    Moreover, $p=q$ if and only if $\sigma$ lies on $M_p$.
\item[(iii)] If $\sigma_1,\ldots,\sigma_l$ are $l$ entire orbits through the random varibales
    $x_1,\ldots,x_l$ respectively such that for some $1\le j_0,\ldots,j_l\le n$,
    $\Omega_{\sigma_k}\subset M_{j_{k-1}}$ and
    $\Omega_{\sigma_k}^*\subset M_{j_{k}}$ for $k=1,\ldots,l$, then $j_0\le j_l$.
    Moreover, $j_0< j_l$ if and only if $\sigma_k$ does not
    lie on $M_D$ with positive probability for some $k$, otherwise $j_0=\cdots=j_l$.
\end{description}
\end{theorem}

\noindent{\bf Proof.} (i) Since
\[
\emptyset=R_0^c\varsubsetneq R_1^c\varsubsetneq\cdots\varsubsetneq
R_n^c=S,
\]
let $\tilde R_i=R_i^c\backslash R_{i-1}^c$. Then
$S=\bigcup_{i=1}^{n}\tilde R_i$ almost surely and $\tilde
R_i=B(A_i)\backslash B(A_{i-1})$. Hence for arbitrary random
variable in $\tilde R_i$, it is attracted by $A_i$ but not by
$A_{i-1}$. For arbitrary random variable $x$ in $S$, choose $n$
random variables $x_1,\ldots,x_n$ such that $x_i\in\tilde R_i$
almost surely and $x(\omega)=x_i(\omega)$ when $\omega\in\Omega_i$,
where $\Omega_i:=\{\omega|x(\omega)\in\tilde R_i(\omega)\}$,
$i=1,\ldots,n$. By the fact $x_i\in\tilde R_i=R_i^c\cap R_{i-1}$ we
know that $x_i$ is attracted by $A_i\cap R_{i-1}=M_i$ almost surely.
Then by Lemma \ref{lem1} we obtain for $\mathbb P$-almost all
$\omega$
\[
\Omega_\sigma(\omega)=\Omega_x(\omega)\subset\bigcup_{i=1}^n\Omega_{x_i}(\omega)
\subset\bigcup_{i=1}^nM_i(\omega)=M_D(\omega).
\]

Since
\[
S=A_0^c\varsupsetneq A_1^c\varsupsetneq\cdots\varsupsetneq
A_n^c=\emptyset,
\]
let $\tilde A_i=A_{i-1}^c\backslash A_{i}^c=B(R_{i-1})\cap A_i$,
$i=1,\ldots,n$. Then $S=\bigcup_{i=1}^{n}\tilde A_i$ almost surely.
By (iv) of Lemma \ref{prop}, for given random variable $x\in\tilde
A_i$, we have $\Omega_\sigma^*\subset R_{i-1}$ almost surely for
{\em any} backward orbit $\sigma$ through $x$. Since $x\in\tilde
A_i\subset A_i$, by the invariance of $A_i$, {\em there exists} a
backward orbit $\sigma$ through $x$ lying on $A_i$ (we can not
guarantee generally that any backward orbit through $x$ must lie on
$A_i$). Hence for this $\sigma$ we have $\Omega_\sigma^*\subset A_i$
almost surely. Therefore, we have obtained that for any random
variable $x\in\tilde A_i$, there exists a backward orbit $\sigma$
through it such that $\Omega_\sigma^*\subset A_i\cap R_{i-1}=M_i$
almost surely. For arbitrary random variable $y\in S$, choose $n$
random variables $y_i$, $i=1,\ldots,n$ such that $y_i\in\tilde A_i$
almost surely and $y(\omega)=y_i(\omega)$ when $\omega\in\Omega_i$,
where $\Omega_i:=\{\omega|y(\omega)\in\tilde A_i(\omega)\}$. By
above argument, for each $i$, there exists a backward orbit
$\sigma_i$ through $y_i$ such that $\Omega_{\sigma_i}^*\subset M_i$
almost surely. ``Attaching" the corresponding parts of these
$\sigma_i$'s together when $y$ lies on $\tilde A_i$, we obtain a
backward orbit $\sigma$ through $y$. By the choice of $\sigma$, we
have
\[
\Omega_\sigma^*\subset\bigcup_{i=1}^n\Omega_{\sigma_i}^*\subset\bigcup_{i=1}^nM_i=M_D
\]
almost surely as desired.

(ii) Since $\Omega_\sigma\subset M_p=A_p\cap R_{p-1}$ almost surely,
we have $x\in A_{p-1}^c$ almost surely. By the fact that
$\Omega_\sigma^*\subset M_q=A_q\cap R_{q-1}$ almost surely, we have
$\sigma$ lying on $A_q$ almost surely by (iii) of Lemma \ref{prop}.
In particular, $\sigma_0=x\in A_q$ almost surely. Hence we have
$x\in A_{p-1}^c\cap A_q$ almost surely. If $q<p$, then $A_q\subset
A_{p-1}$, hence $A_q\cap A_{p-1}^c=\emptyset$ almost surely, a
contradiction.

If $\sigma$ lies on $M_p$, then we have
$\Omega_\sigma,\Omega_\sigma^*\subset M_p$ almost surely by the fact
that $M_p$ is an invariant random compact set. That is, we must have
$p=q$. Conversely, if $p=q$, the fact $\Omega_\sigma=\Omega_x\subset
M_p=A_p\cap R_{p-1}$ implies $x\in R_{p-1}$ almost surely. It
follows that $\sigma$ lies on $R_{p-1}$ since any backward orbit
through a random variable in $R_{p-1}$ must lie on it.
$\Omega_\sigma^*\subset M_p=A_p\cap R_{p-1}$ implies that $\sigma$
lies on $A_{p}$ by (iii) of Lemma \ref{prop}. So we have obtained
that $\sigma$ lies on $A_p\cap R_{p-1}=M_p$ almost surely.

(iii) follows from (ii) immediately. \hfill$\Box$

\begin{remark}\rm
In (i) of Theorem \ref{Mor}, we obtain that, for given random
variable $x$, there exists an entire orbit through $x$ satisfying
$\Omega_\sigma\subset M_D$ and $\Omega_\sigma^*\subset M_D$ almost
surely. While in the deterministic case, any entire orbit has this
property, see \cite{Con} for the flow case and \cite{Ryb} for the
semiflow case. But their methods are not applicable here. We are not
sure whether the similar result holds for random semiflow.
\end{remark}

In Theorem \ref{Mor}, we give the positive answer to an open problem
put forward by Caraballo and Langa \cite{CL}, which confirms our
belief that we may understand theoretically the structure of global
random attractor in a way as we do in the deterministic case. But,
as mentioned in introduction, the theory of dynamics on the global
random attractor is so immature that we know very few examples for
which the structure of global random attractors is well understood.
Hence to understand them well, more concrete examples are needed at
the first step. The global random attractor of Chafee-Infante
reaction-diffusion equation perturbed by Stratonovich multiplicative
noise may be the best understood model in infinite dimensional case,
so we use it to illustrate our results. For this model in
deterministic case, it is well studied, see Hale \cite{Hale} and
Henry \cite{Hen} for instance.

\begin{example}\rm
Consider the Chafee-Infante reaction-diffusion equation perturbed by
Stratonovich multiplicative noise
\begin{equation}\label{eq}
{\rm d} u=(\Delta u+\beta u-u^3){\rm d} t+\delta u\circ{\rm d}
W,~~x\in[0,\pi]~{\rm with}~u(t,0)=u(t,\pi)=0.
\end{equation}
This equation is well studied in \cite{CLR} and very recently in
\cite{CCLR} and \cite{WD}. Here we take it as an example to
illustrate our results, for detailed analysis see \cite{CLR}. Denote
\begin{eqnarray*}
&D=[0,\pi],~~H=L^2(D),\\
&\mathcal K^+=\{u\in H| u(x)\ge 0 ~{\rm almost~ everywhere}\},\\
&\mathcal K^-=\{u\in H| u(x)\le 0 ~{\rm almost~ everywhere}\}.
\end{eqnarray*}
Assume that $\varphi$ is the RDS generated by (\ref{eq}) on $H$,
then $\varphi$ is order preserving on $H$ \cite{Kot}, i.e. if
$u_0\ge v_0$ almost everywhere, then
\[
\varphi(t,\omega)u_0\ge \varphi(t,\omega)v_0.
\]
In particular, since $\{0\}$ is a random fixed point of $\varphi$,
$\mathcal K^\pm$ are invariant cones of $\varphi$. We can show that
in each of these two cones, there exists a random compact absorbing
set. Hence it follows that $\varphi$ has nontrivial random
attractors $\mathcal A^+$, $\mathcal A^-$, respectively, in each of
these two cones. Furthermore, there exist positive and negative
random fixed points $\pm a(\omega)\in\mathcal A^\pm(\omega)$ such
that
\begin{eqnarray*}
& 0\le u\le a(\omega)&{\rm for ~all ~}u\in\mathcal A^+(\omega),\\
-&a(\omega)\le u\le 0\quad &{\rm for ~all ~}u\in\mathcal
A^-(\omega).
\end{eqnarray*}
It has been shown in \cite{CLR} that (\ref{eq}) has global random
attractor $\mathcal A$ and the lower bound on the dimension of
$\mathcal A$ has also been obtained, but the exact structure of
$\mathcal A$ has not been obtained in \cite{CLR}. The authors
conjecture in \cite{CLR} that, for $\lambda_1<\beta<\lambda_2$
($\lambda_1$, $\lambda_2$ are the first two eigenvalues of the
negative Laplacian), the global random attractor
\[
\mathcal A(\omega)=\mathcal A^+(\omega)\cup\mathcal A^-(\omega),
\]
with $\mathcal A^\pm$ consisting of a one-dimensional manifold
joining the origin to $\pm a(\omega)$. Moreover, since orbits near
$0$ in $\mathcal K^+$ move away from the origin (i.e. $\{0\}$ is
unstable), it is also conjectured in \cite{CLR} that $a(\omega)$ is
attracting in $\mathcal K^+$ (symmetrically, $-a(\omega)$ is
attracting in $\mathcal K^-$). Very recently, Wang and Duan
\cite{WD} confirm these two conjectures. Indeed, they show that the
semiflow on the global random attractor of (\ref{eq}) is
topologically equivalent to a well studied one-dimensional
stochastic ODE (see \cite{AB,CIS} for the study of this stochastic
ODE), see \cite{WD} for details. Now we can use our abstract results
here. Choose $\pm b(\omega)$ such that $0<b(\omega)<a(\omega)$ and
$-a(\omega)<-b(\omega)<0$, then by the order preserving property of
$\varphi$ and the fact that $\pm a(\omega)$ are attracting in
$\mathcal K^\pm$, respectively, we have $[-a(\omega),-b(\omega))$
and $(b(\omega),a(\omega)]$ being forward invariant and
$\Omega_{[-a,-b)}(\omega)=\{-a(\omega)\}$,
$\Omega_{(b,a]}(\omega)=\{a(\omega)\}$. Therefore, it is clear that
$\emptyset$, $\{a(\omega)\}$, $\{-a(\omega)\}$,
$\{a(\omega),-a(\omega)\}$ and $\mathcal A$ are all local random
attractors on the global random attractor $\mathcal A$. If we set
\[
A_0=\emptyset,~A_1(\omega)=\{a(\omega)\},~A_2(\omega)=\{a(\omega),-a(\omega)\},~A_3(\omega)=\mathcal
A(\omega),
\]
then the corresponding repellers are
\[
R_0(\omega)=\mathcal
A(\omega),~R_1(\omega)=[-a(\omega),0],~R_2(\omega)=\{0\},~R_3=\emptyset,
\]
respectively. So $D=\{M_1,M_2,M_3\}$ with
\begin{equation}\label{Mo}
M_1=\{a(\omega)\},~~M_2=\{-a(\omega)\},~~M_3=\{0\}
\end{equation}
is a Morse decomposition of $\mathcal A$, just consisting of random
fixed points of $\varphi$. If we reset $A_0(\omega)=\emptyset$,
$A_1=\{-a(\omega)\}$, $A_2(\omega)=\{a(\omega),-a(\omega)\}$,
$A_3(\omega)= \mathcal A(\omega)$, then
$M_1(\omega)=\{-a(\omega)\}$, $M_2(\omega)=\{a(\omega)\}$,
$M_3(\omega)=\{0\}$ are the corresponding Morse sets. These two are
the finest Morse decompositions of $\mathcal A$, consisting of just
random fixed points. The Morse decomposition can be coarsen: if we
set $A_0=\emptyset$, $A_1(\omega)=\{-a(\omega),a(\omega)\}$,
$A_2=\mathcal A$, then $R_0=\mathcal A$, $R_1(\omega)=\{0\}$,
$R_2=\emptyset$. Hence the corresponding Morse sets are
$M_1=\{-a(\omega),a(\omega)\}$, $M_2=\{0\}$.

By Remark \ref{con} we know that the random chain recurrent set on
$\mathcal A$
\[
\mathcal{CR}_\varphi(\omega)=\bigcap[A(\omega)\cup
R(\omega)]=\{-a(\omega),0,a(\omega)\}
\]
for $\mathbb P$-almost all $\omega$, where the intersection is taken
over all local random attractors on $\mathcal A$, i.e. $\emptyset$,
$\{a(\omega)\}$, $\{-a(\omega)\}$, $\{a(\omega),-a(\omega)\}$ and
$\mathcal A$. It is clear that the random chain recurrent set on
$\mathcal A$ equals the union of Morse sets in (\ref{Mo}).
\end{example}

\vskip8mm

\noindent{\bf Acknowledgements}\\
The author expresses his sincere thanks to Professor Yong Li for his
instruction and many invaluable suggestions. The author is grateful
to Professor Jack Hale and Professor Yingfei Yi for their proposal
on this problem. The author sincerely thanks the anonymous referees
for their careful reading the manuscript and invaluable comments
which are very helpful to improve the paper.

\end{document}